\documentclass{amsart}%
\usepackage{graphicx, amssymb, amsfonts, amsmath}
\usepackage{amsmath}
\usepackage{amsfonts}
\usepackage{amssymb}
\usepackage{graphicx}%
\providecommand{\U}[1]{\protect\rule{.1in}{.1in}}
\providecommand{\U}[1]{\protect\rule{.1in}{.1in}}
\newtheorem{theorem}{Theorem}
\theoremstyle{plain}

\newtheorem{proposition}{Proposition}

\numberwithin{equation}{section}
\begin{document}
\title{A Topologically Induced 2-in/2-out Operation on Loop Cohomology }
\author{Ronald Umble$^{1}$}
\address{Department of Mathematics\\
Millersville University of Pennsylvania\\
Millersville, PA. 17551}
\email{ron.umble@millersville.edu}
\thanks{$^{1}$ This research funded in part by a Millersville University faculty
research grant.}
\date{September 15, 2011; revised March 21, 2012}

\begin{abstract}
We apply the Transfer Algorithm introduced in \cite{SU5} to transfer an
$A_{\infty}$-algebra structure that cannot be computed using the classical
Basic Perturbation Lemma. We construct a space $X$ whose topology induces a
nontrivial 2-in/2-out operation $\omega_{2}^{2}$ on loop cohomology $H^{\ast
}\left(  \Omega X;\mathbb{Z}_{2}\right)  $.

\end{abstract}
\maketitle

\section{Introduction}

In \cite{SU4} and \cite{SU5}, S. Saneblidze and this author defined the
notions of a matrad and a relative matrad, and constructed the related
families of polytopes known as biassociahedra $KK=\left\{  KK_{n,m}%
=KK_{m,n}\right\}  $ and bimultiplihedra $JJ=\left\{  JJ_{n,m}=JJ_{m,n}%
\right\}  $ of which $KK_{1,n}$ is the associahedron $K_{n}$ and $JJ_{1,n}$ is
the multiplihedron $J_{n}.$ Cells of $KK$ and $JJ$ are identified with certain
fraction product monomials, for example,\smallskip%

\[
\raisebox{-0.2in}{\includegraphics[
trim=0.000000in -0.113930in 0.000000in 0.113930in,
height=.45in,
width=1.3465in
]{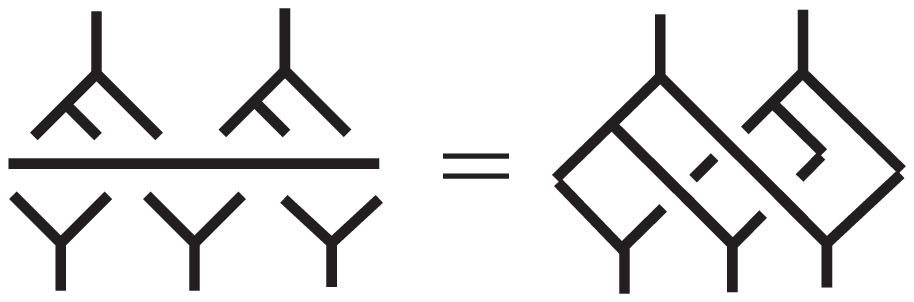}} \text{ \ }\leftrightarrow\text{ \ a vertex of }KK_{2,3}.
\]
\smallskip

\noindent In fact, $JJ_{m,n}$ is a subdivision of $KK_{m,n}\times I\ $with
$\partial JJ_{m,n}$ containing the cells $KK_{n,m}\times0$ and $KK_{n,m}%
\times1$.

Let $R$ be a commutative ring with unity. The free matrad $\mathcal{H}%
_{\infty}$ is represented by the DG $R$-module (DGM) of cellular chains
$C_{\ast}\left(  KK\right)  $ by associating the top dimensional cell of
$KK_{n,m}$ with the matrad generator $\theta_{m}^{n}\in\mathcal{H}_{\infty}%
:$\smallskip%

\[%
\begin{array}
[c]{cc}%
n\text{ outputs} & \\
\raisebox{-0.2361in}{\includegraphics[
height=0.5474in,
width=0.5561in
]{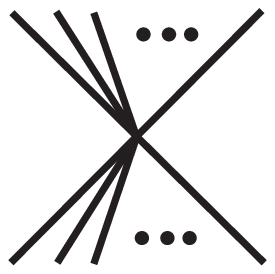}} & \leftrightarrow\text{ \ }\theta_{m}^{n}\text{\ }.\\
m\text{ inputs} &
\end{array}
\]
\smallskip

\noindent An $A_{\infty}$\emph{-bialgebra} is a DGM $(A,d)$ together with a
family of multilinear operations $\omega=\{\omega_{m}^{n}\in Hom^{m+n-3}%
(A^{\otimes m},A^{\otimes n})\mid mn\neq1\}$ and a map of matrads
${\mathcal{H}}_{\infty}\rightarrow\mathcal{E}nd_{TA}$ such that $\theta
_{m}^{n}\mapsto\omega_{m}^{n}$, i.e., $\left(  A,\omega\right)  $ is an
algebra over $\mathcal{H}_{\infty}.$ Note that we recover the operadic
structure of $A_{\infty}$-(co)algebras by setting $m=1$ or $n=1$.

Similarly, the\ free relative matrad $\mathcal{JJ}_{\infty}$ is represented by
the DGM$\ $of cellular chains $C_{\ast}\left(  JJ\right)  $ by associating the
top dimensional cell of $JJ_{m,n}$ with the relative matrad generator
$\mathfrak{f}_{m}^{n}\in\mathcal{JJ}_{\infty}:$%

\[%
\begin{array}
[c]{cc}%
n\text{ outputs} & \\
\raisebox{-0.2361in}{\includegraphics[
height=0.5474in,
width=0.5561in
]{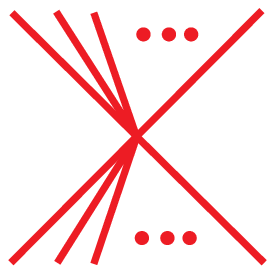}} & \leftrightarrow\text{ \ }\mathfrak{f}_{m}^{n}%
\text{\ }.\\
m\text{ inputs} &
\end{array}
\]
\smallskip

\noindent Let $(A,\omega_{A})$ and $(B,\omega_{B})$ be $A_{\infty}%
$-bialgebras. A \emph{morphism} $G$ \emph{from} $A$ \emph{to} $B,$ denoted by
$G:A\Rightarrow B,$ is a family of multilinear maps $G=\{g_{m}^{n}\in
Hom^{m+n-2}(A^{\otimes m},B^{\otimes n})\}$ together with a map of relative
matrads $\mathcal{J}\mathcal{J}_{\infty}\rightarrow Hom\left(  TA,TB\right)  $
such that $\mathfrak{f}_{m}^{n}\mapsto g_{m}^{n}$, i.e., $G$ is an
$\mathcal{H}_{\infty}$-bimodule. The elements $\theta_{m}^{n}\left(
\mathfrak{f}_{1}^{1}\right)  ^{\otimes m}$ and $\left(  \mathfrak{f}_{1}%
^{1}\right)  ^{\otimes n}\theta_{m}^{n}$ of $\mathcal{J}\mathcal{J}_{\infty}$
are associated with the codimension 1 cells $KK_{n,m}\times0$ and
$KK_{n,m}\times1$ of $JJ_{m,n},$ respectively,

\begin{center}
$%
\begin{array}
[c]{ccccccc}%
n &  &  &  & n &  & \\
\raisebox{-0.2707in}{\includegraphics[
height=0.6106in,
width=0.5025in
]{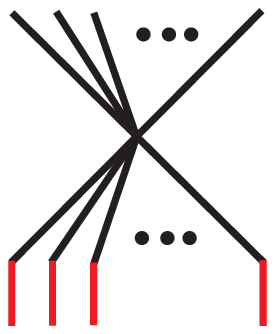}} & \leftrightarrow & \theta_{m}^{n}\left(  \mathfrak{f}%
_{1}^{1}\right)  ^{\otimes m}\text{\ }; &  &
\raisebox{-0.2569in}{\includegraphics[
height=0.6036in,
width=0.5025in
]{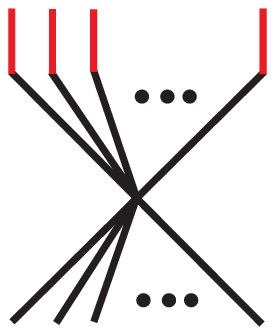}} & \leftrightarrow & \left(  \mathfrak{f}_{1}^{1}\right)
^{\otimes n}\theta_{m}^{n}\text{\ },\\
m &  &  &  & m &  &
\end{array}
$
\end{center}

\noindent and the aforementioned map of relative matrads sends $\theta_{m}%
^{n}\left(  \mathfrak{f}_{1}^{1}\right)  ^{\otimes m}\mapsto\omega_{m}%
^{n}g^{\otimes m}$ and $\left(  \mathfrak{f}_{1}^{1}\right)  ^{\otimes
n}\theta_{m}^{n}\mapsto g^{\otimes n}\omega_{m}^{n} $. Again, we recover the
structure of an $A_{\infty}$-(co)algebra morphism by setting $m=1$ or $n=1$. A
morphism $\Phi=\{\phi_{m}^{n}\}_{m,n\geq1}:A\Rightarrow B$ is an
\emph{isomorphism} if $\phi_{1}^{1}{\ }$is an isomorphism of underlying modules.

The paper is organized as follows: In section 2 we review the \emph{Transfer
Algorithm }introduced in \cite{SU5} and apply it to transfer an $A_{\infty}%
$-algebra structure that cannot be computed using the classical Basic
Perturbation Lemma. In Section 3 we construct a space $X$ whose topology
induces a nontrivial 2-in/2-out operation $\omega_{2}^{2}$ on loop cohomology
$H^{\ast}\left(  \Omega X;\mathbb{Z}_{2}\right)  $.

\section{Transfer of $A_{\infty}$-Structure}

If $A$ is a free DGM, $B$ is an $A_{\infty}$-algebra, and $g:A\rightarrow B$
is a homology isomorphism with a right-homotopy inverse, the Basic
Perturbation Lemma (BPL) transfers the $A_{\infty}$-algebra structure from $B
$ to $A$ (see \cite{Hueb-Kade}, \cite{Markl1}, for example). When $B$ is an
$A_{\infty}$-bialgebra, Theorem \ref{transfer} generalizes the BPL in two directions:

\begin{enumerate}
\item The $A_{\infty}$-bialgebra structure on $B$ transfers to an $A_{\infty}
$-bialgebra structure on $A.$

\item The transfer algorithm requires neither freeness in $A$ nor the
existence of a right-homotopy inverse of $g$.
\end{enumerate}

Given DGMs $\left(  A,d_{A}\right)  $ and $\left(  B,d_{B}\right)  ,$ let
$\nabla$ be the induced differential on $U_{A,B}=Hom\left(  TA,TB\right)  $,
i.e., for $f\in U_{A,B}$ define $\nabla f=d_{B} f-\left(  -1\right)
^{\left\vert f\right\vert }f d_{A},$ where $d_{A}$ and $d_{B}$ denote the free
linear extensions of $d_{A}$ and $d_{B}$. A chain map $g:A\rightarrow B$
induces a cochain map $\tilde{g}:\mathcal{E}nd_{TA}\rightarrow U_{A,B}$
defined on $u\in Hom\left(  A^{\otimes m},A^{\otimes n}\right)  $ by
$\tilde{g}\left(  u\right)  =g^{\otimes n}u.$ If $g$ is a homology
isomorphism, so is $\tilde{g} $ provided condition (i) or (ii) in the
following proposition is satisfied (the proof is left to the reader):

\begin{proposition}
\label{homology-iso} Let $\left(  A,d_{A}\right)  $ and $\left(
B,d_{B}\right)  $ be DGMs, and let $g:A\rightarrow B$ be a chain map that is
also a homology isomorphism. Then $\tilde{g}:\mathcal{E}nd_{TA}\rightarrow
U_{A,B}$ is a homology isomorphism if either of the following conditions holds:

\begin{enumerate}
\item[\textit{(i)}] $A$ is free as an $R$-module.

\item[\textit{(ii)}] For each $n\geq1,$ there is a DGM\ $X\left(  n\right)  $
and a splitting $B^{\otimes n}=A^{\otimes n}\oplus X(n)$ as a chain complex
such that $H^{\ast}Hom\left(  A^{\otimes k},X\left(  n\right)  \right)  =0$
for all $k\geq1.$
\end{enumerate}
\end{proposition}

Thus there is the following generalization of the BPL:

\begin{theorem}
[\textbf{The} \textbf{Transfer}]\label{transfer}Let $\left(  A,d_{A}\right)  $
be a DGM, let $(B,d_{B},\omega_{B})$ be an $A_{\infty}$-bialgebra, and let
$g:A\rightarrow B $ be a chain map and a homology isomorphism. If $\tilde
{g}:\mathcal{E}nd_{TA}\rightarrow U_{A,B}$ is a homology isomorphism, then

\begin{enumerate}
\item[\textit{(i)}] (Existence) $g$ induces an $A_{\infty}$-bialgebra
structure $\omega_{A}=\{\omega_{A}^{n,m}\}$ on $A$ and extends to a map
$G=\{g_{m}^{n}\mid g_{1}^{1}=g\}:A\Rightarrow B$\ of $A_{\infty}$-bialgebras.

\item[\textit{(ii)}] (Uniqueness) $\left(  \omega_{A},G\right)  $ is unique up
to isomorphism, i.e., if $\left(  \omega_{A},G\right)  $ and $\left(
\bar{\omega}_{A},\bar{G}\right)  $ are induced by chain homotopic maps $g$ and
$\bar{g}$, there is an isomorphism $\Phi:\left(  A,\bar{\omega}_{A}\right)
\Rightarrow\left(  A,\omega_{A}\right)  $ and a chain homotopy $T:\bar
{G}\simeq G\circ\Phi.$
\end{enumerate}
\end{theorem}

The proof of Theorem \ref{transfer}, which appears in \cite{SU5}, suggests the
following general \emph{Transfer Algorithm}:%

\[%
\begin{tabular}
[c]{l}%
\textbf{The} \textbf{Transfer Algorithm }%
\end{tabular}
\ \
\]

\noindent\textbf{Initial data }

\begin{itemize}
\item A DGM $\left(  A,d_{A}\right)  $

\item An $A_{\infty}$-bialgebra $\left(  B,d_{B},\omega_{B}\right)  $ and a
map of matrads $\alpha_{B}:C_{\ast}(KK)\rightarrow\mathcal{E}nd_{TB}$ sending
$\theta_{m}^{n}\mapsto\omega_{B}^{n,m}$

\item A chain map/homology isomorphism $g:A\rightarrow B$ such that $\tilde
{g}$ is a homology isomorphism$\smallskip$
\end{itemize}

\noindent\textbf{Objectives}

\begin{itemize}
\item Define operations $\omega_{A}^{n,m}:A^{\otimes m}\rightarrow A^{\otimes
n}$ for all $m,n,$ $mn\neq1$

\item Construct a map of matrads $\alpha_{A}:C_{\ast}(KK)\rightarrow
\mathcal{E}nd_{TA}$ sending $\theta_{m}^{n}\mapsto\omega_{A}^{n,m}$

\item Construct a map of $A_{\infty}$-bialgebras $G=\left\{  g_{m}^{n}\mid
g_{1}^{1}=g\right\}  :A\Rightarrow B\smallskip$
\end{itemize}

\noindent\textbf{Initialization}

\begin{enumerate}
\item[1.] Define $\beta:C_{0}\left(  JJ_{1,1}\right)  \rightarrow Hom\left(
A,B\right)  $ by $\mathfrak{f}_{1}^{1}=
\raisebox{-0.0623in}{\includegraphics[clip=true, viewport=.0in .0in .1in .1in,
trim=0.0000000in 0.0000000in 0.0000000in 0.302082in,
height=0.2067in,
width=0.2551in
]{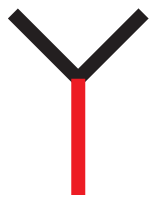}} \ \mapsto g\vspace{.05in}$

\item[2.] Define $\beta$ on the vertex
$\raisebox{-0.0553in}{\includegraphics[
height=0.18in,
width=0.15in
]{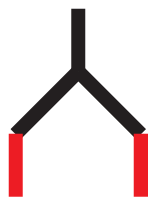}} $ of $JJ_{1,2}$ by $\theta_{2}^{1}\left(  \mathfrak{f}%
_{1}^{1}\otimes\mathfrak{f}_{1}^{1}\right)  \ \mapsto\omega_{B}^{1,2}\left(
g\otimes g\right)  \vspace{.05in}$

\item[3.] Define $\beta$ on the vertex
$\raisebox{-0.0484in}{\includegraphics[
height=0.18in,
width=0.15in
]{beta1-2-rb.eps}} $ of $JJ_{2,1}$ by $\theta_{1}^{2}\mathfrak{f}_{1}%
^{1}\ \mapsto\omega_{B}^{2,1}g\vspace{.05in}$

\item[4.] Consider the $\nabla$-cocycle $\omega_{B}^{1,2}\left(  g\otimes
g\right)  $

\begin{itemize}
\item Choose a cocycle $\omega_{A}^{1,2}\in\mathcal{E}nd_{TA}$ such that
$\tilde{g}_{\ast}[\omega_{A}^{1,2}]=[\omega_{B}^{1,2}\left(  g\otimes
g\right)  ]$

\item Define $\alpha_{A}:C_{0}\left(  KK_{1,2}\right)  \rightarrow Hom\left(
A^{\otimes2},A\right)  $ by $\theta_{2}^{1}=
\raisebox{-0.0623in}{\includegraphics[
height=0.19in,
width=0.20in
]{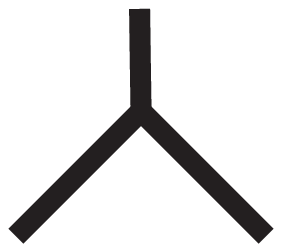}} \ \mapsto\omega_{A}^{1,2}$

\item Define $\alpha_{A}:C_{0}\left(  \partial KK_{1,3}\right)  \rightarrow
Hom\left(  A^{\otimes3},A\right)  $ by%
\[
\raisebox{-0.083in}{\includegraphics[
height=0.2525in,
width=0.2586in
]{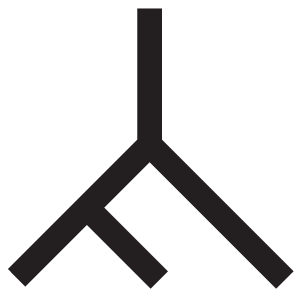}} \ \mapsto\omega_{A}^{1,2}\left(  \omega_{A}^{1,2}%
\otimes\mathbf{1}\right)  \ \text{and}%
\ \raisebox{-0.0969in}{\includegraphics[
height=0.2525in,
width=0.2586in
]{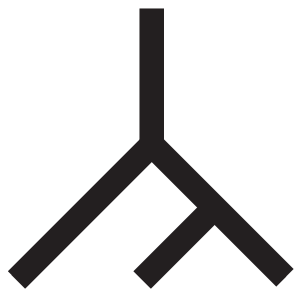}} \ \mapsto\omega_{A}^{1,2}\left(  \mathbf{1}\otimes\omega
_{A}^{1,2}\right)
\]

\item Extend $\beta$ to the vertex $\raisebox{-0.05in}{\includegraphics[
height=0.18in,
width=0.15in
]{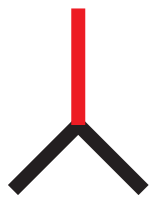}} \subset JJ_{1,2}$ via $\mathfrak{f}_{1}^{1}\theta_{2}%
^{1}\ \mapsto g\omega_{A}^{1,2}$
\end{itemize}

\item[5.] Dually, consider the $\nabla$-cocycle $\omega_{B}^{2,1}g$

\begin{itemize}
\item Choose a cocycle $\omega_{A}^{2,1}\in\mathcal{E}nd_{TA}$ such that
$\tilde{g}_{\ast}[\omega_{A}^{2,1}]=[\omega_{B}^{2,1}g]$

\item Define $\alpha_{A}:C_{0}\left(  KK_{2,1}\right)  \rightarrow Hom\left(
A,A^{\otimes2}\right)  $ by $\theta_{1}^{2}=
\raisebox{-0.0692in}{\includegraphics[
height=0.19in,
width=0.20in
]{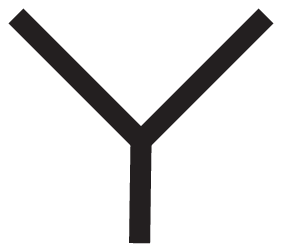}} \ \mapsto\omega_{A}^{2,1}$

\item Define $\alpha_{A}:C_{0}\left(  \partial KK_{3,1}\right)  \rightarrow
Hom\left(  A,A^{\otimes3}\right)  $ by%
\[
\ \raisebox{-0.083in}{\includegraphics[
height=0.256in,
width=0.2586in
]{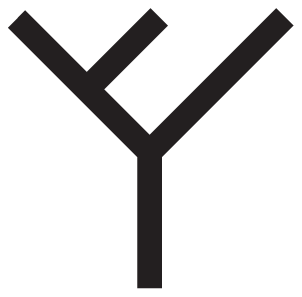}} \ \mapsto\left(  \omega_{A}^{2,1}\otimes\mathbf{1}\right)
\omega_{A}^{2,1}\ \text{and}\ \raisebox{-0.0899in}{\includegraphics[
height=0.256in,
width=0.2586in
]{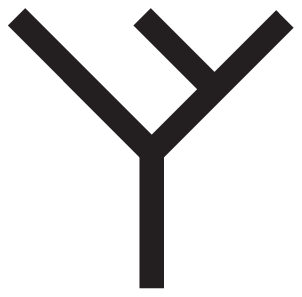}} \ \mapsto\left(  \mathbf{1}\otimes\omega_{A}^{2,1}\right)
\omega_{A}^{2,1}%
\]

\item Extend $\beta$ to the vertex $\raisebox{-0.05in}{\includegraphics[
height=0.18in,
width=0.15in
]{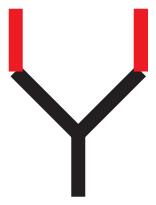}} \subset JJ_{1,2}$ via $\left(  \mathfrak{f}_{1}^{1}%
\otimes\mathfrak{f}_{1}^{1}\right)  \theta_{1}^{2}\ \mapsto\left(  g\otimes
g\right)  \omega_{A}^{2,1}$
\end{itemize}

\item[6.] Define $\alpha_{A}:C_{0}\left(  \partial KK_{2,2}\right)
\rightarrow Hom\left(  A^{\otimes2},A^{\otimes2}\right)  $ by%
\[
\hspace{.4in} \raisebox{-0.1323in}{\includegraphics[
height=0.3425in,
width=0.4272in
]{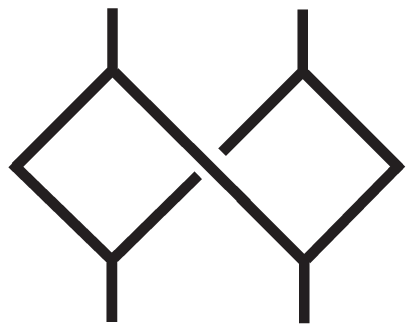}} \ \mapsto\left(  \omega_{A}^{1,2}\otimes\omega
_{A}^{1,2}\right)  \sigma_{2,2}\left(  \omega_{A}^{2,1}\otimes\omega_{A}%
^{2,1}\right)  \text{\ and\ } \raisebox{-0.1038in}{\includegraphics[
height=0.2741in,
width=0.1903in
]{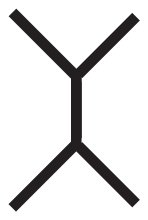}} \ \mapsto\omega_{A}^{2,1}\omega_{A}^{1,2},\text{\ where}%
\]
$\sigma_{p,q}:(A^{\otimes p})^{\otimes q} \overset{\approx}{\to} (A^{\otimes
q})^{\otimes p}$ is the canonical permutation of tensor factors

\item[7.] Note that $\left[  \omega_{B}^{1,2}\left(  g\otimes g\right)
-g\omega_{A}^{1,2}\right]  =0$

\begin{itemize}
\item Choose a cochain $g_{2}^{1}$ such that $\nabla g_{2}^{1}=\omega
_{B}^{1,2}\left(  g\otimes g\right)  -g\omega_{A}^{1,2}$

\item Define $\beta:C_{1}\left(  JJ_{1,2}\right)  \rightarrow Hom\left(
A^{\otimes2},B\right)  $ by $\mathfrak{f}_{2}^{1}=$
$\raisebox{-0.0484in}{\includegraphics[
height=0.1781in,
width=0.1384in
]{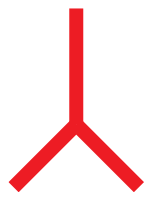}} \ \mapsto g_{2}^{1}$

\item Define $\beta$ on a monomial in $C_{\ast}\left(  \partial JJ_{1,3}%
\smallsetminus\text{int}KK_{1,3}\times1\right)  $ to be the corresponding
composition:$\medskip$
\[
\raisebox{0.0899in}{\includegraphics[
trim=0.000000in -0.113930in 0.000000in 0.113930in,
height=1.5826in,
width=1.3465in
]{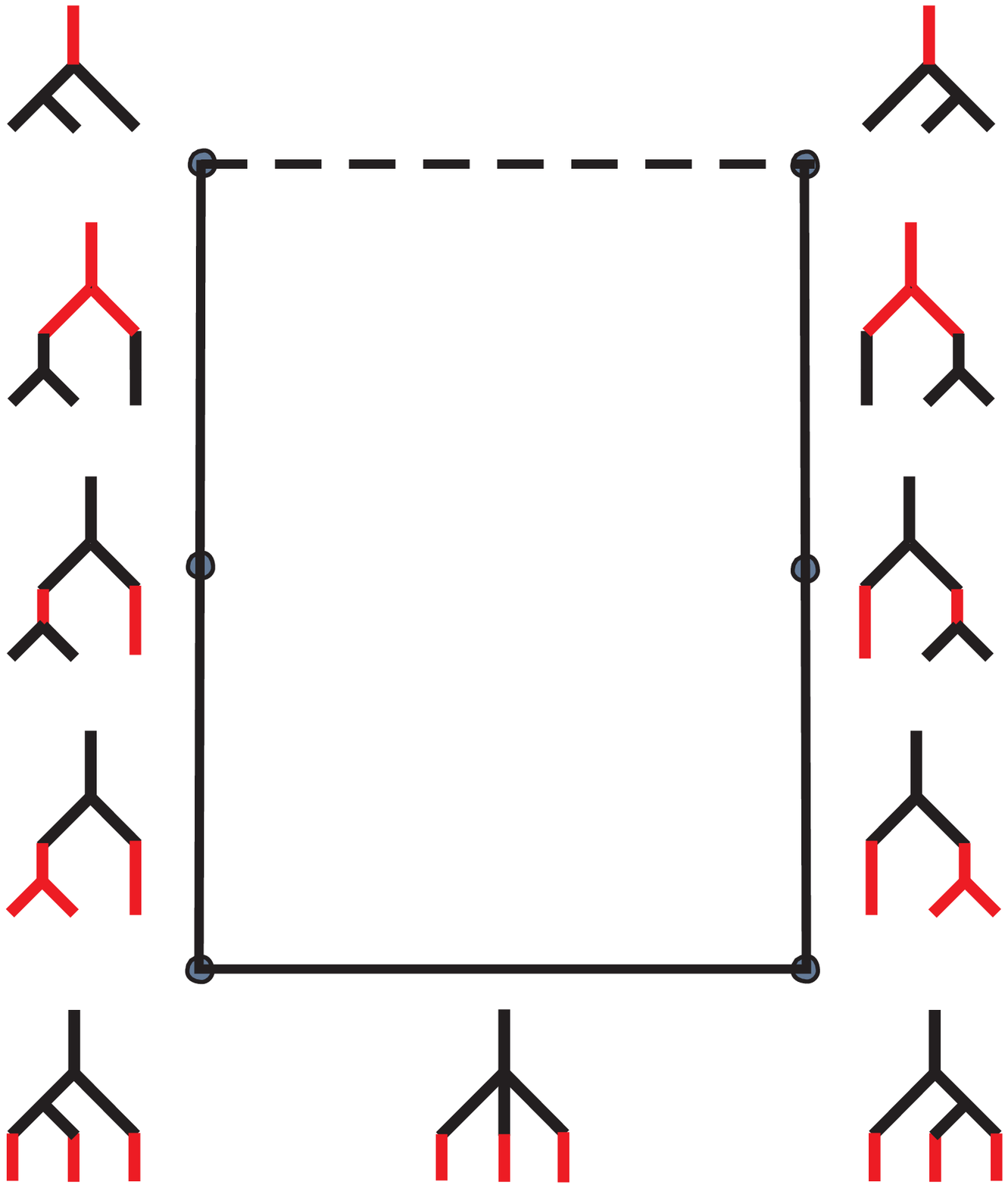}}
\]

\end{itemize}

\item[8.] Dually, note that $\left[  \omega_{B}^{2,1}g-\left(  g\otimes
g\right)  \omega_{A}^{2,1}\right]  =0$

\begin{itemize}
\item Choose a cochain $g_{1}^{2}\in U_{A,B}$ such that $\nabla g_{1}%
^{2}=\omega_{B}^{2,1}g-\left(  g\otimes g\right)  \omega_{A}^{2,1}$

\item Define $\beta:C_{1}\left(  JJ_{2,1}\right)  \rightarrow Hom\left(
A,B^{\otimes2}\right)  $ by $\mathfrak{f}_{1}^{2}=$
$\raisebox{-0.0553in}{\includegraphics[
height=0.1781in,
width=0.1384in
]{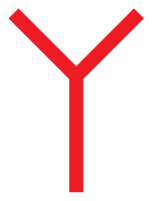}} \ \mapsto g_{1}^{2}$

\item Define $\beta$ on a monomial in $C_{\ast}\left(  \partial JJ_{3,1}%
\smallsetminus\text{int}KK_{3,1}\times1\right)  $ to be the corresponding
composition:$\medskip$%
\[
\raisebox{0.1038in}{\includegraphics[
height=1.5826in,
width=1.3465in
]{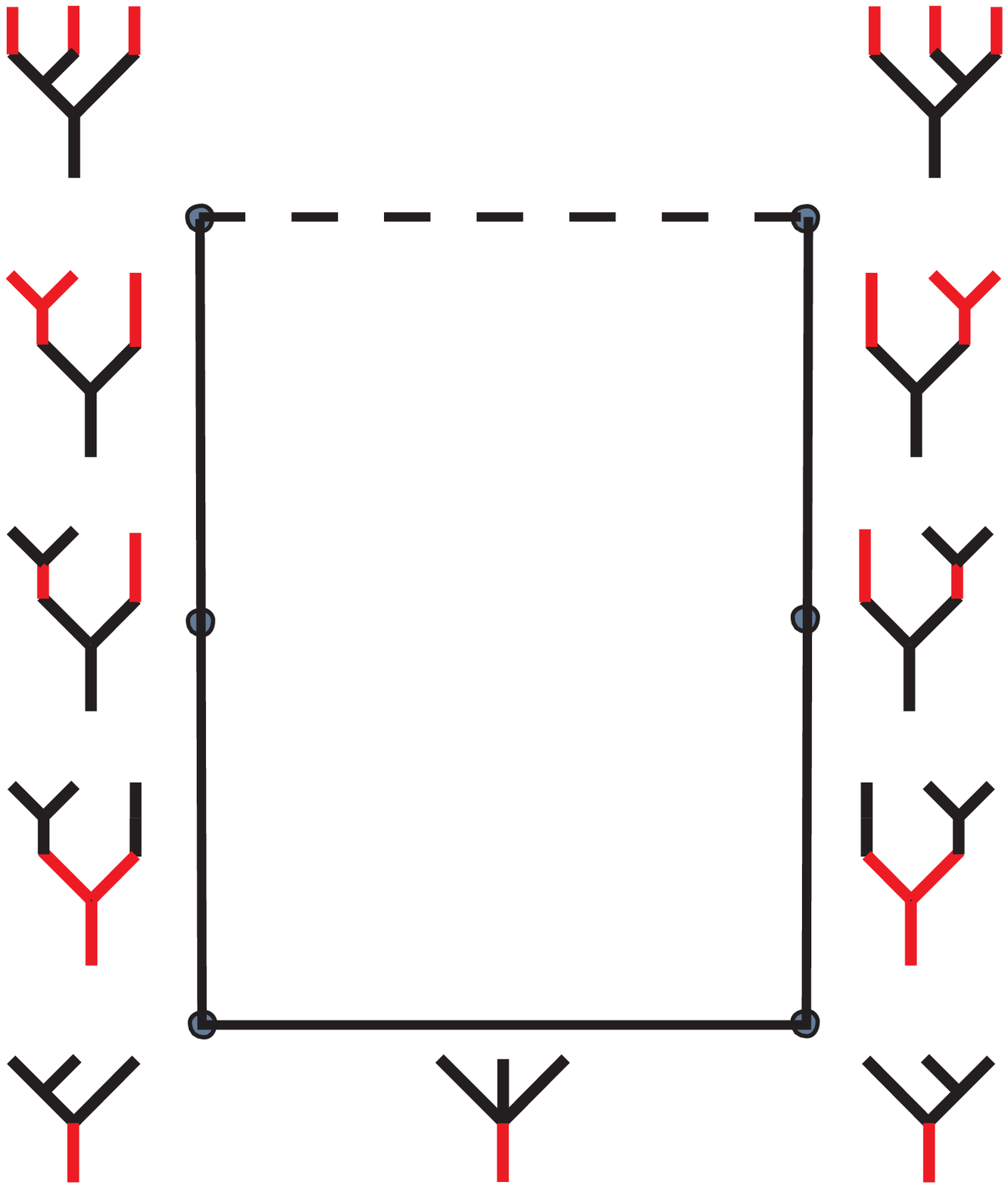}}
\]

\end{itemize}

\item[9.] Define $\beta$ on a monomial in $C_{\ast}\left(  \partial
JJ_{2,2}\smallsetminus\text{int}KK_{2,2}\times1\right)  $ to be the
corresponding fraction product:
\end{enumerate}

\[
\text{ \raisebox{0.0138in}{\includegraphics[
height=1.7469in,
width=1.8299in
]{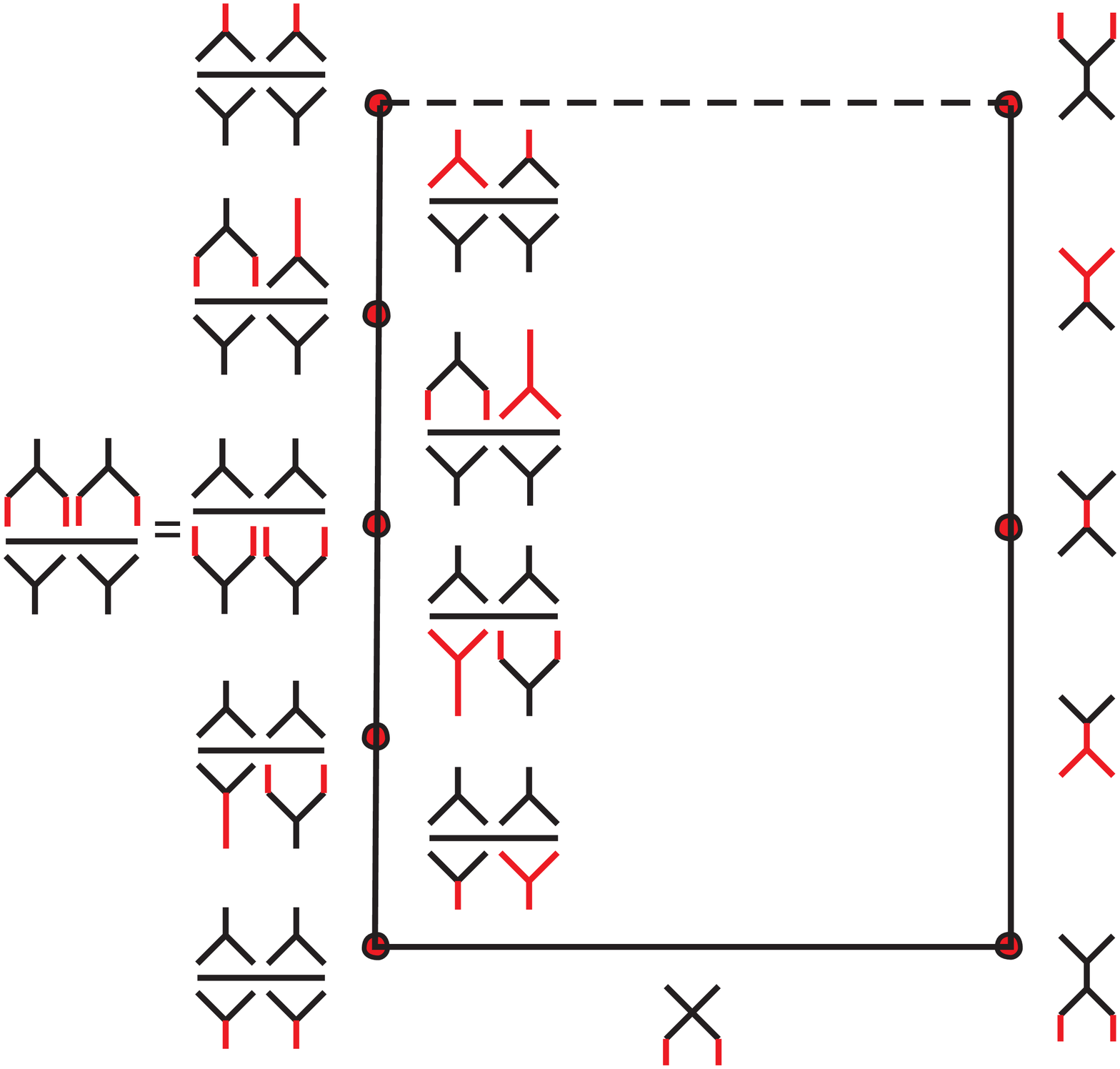}} }%
\]

\noindent\textbf{Induction hypothesis}

Given $m+n\geq4$, assume that for $i+j<m+n,$ $ij\neq1,$ there exists a map

\begin{itemize}
\item $\alpha_{A}:C_{\ast}\left(  KK_{j,i}\right)  \rightarrow Hom\left(
A^{\otimes i},A^{\otimes j}\right)  $ of matrads sending $\theta_{i}%
^{j}\mapsto\omega_{A}^{j,i}$

\item $\beta:C_{\ast}\left(  JJ_{j,i}\right)  \rightarrow Hom\left(
A^{\otimes i},B^{\otimes j}\right)  $ of relative matrads sending
$\mathfrak{f}_{i}^{j}\mapsto g_{i}^{j}$\smallskip
\end{itemize}

\noindent\textbf{Induction objectives}

\begin{itemize}
\item Define $\alpha_{A}$ on the generator $\theta_{m}^{n} \in C_{m+n-3}%
\left(  KK_{n,m}\right)  $

\item Define $\beta$ on the monomial $\left(  \mathfrak{f}_{1}^{1}\right)
^{\otimes n}\theta_{m}^{n}\in C_{m+n-3}\left(  JJ_{n,m}\right)  $

\item Define $\beta$ on the generator $\mathfrak{f}_{m}^{n}\in C_{m+n-2}%
\left(  JJ_{n,m}\right)  $
\end{itemize}

\smallskip

\noindent\textbf{Induction}

For each $i+j=m+n,$ $ij\neq1$

\begin{enumerate}
\item[1.] Define $\alpha_{A}$ on each monomial in $C_{\ast}\left(  \partial
KK_{n,m}\right)  $ to be its corresponding fraction product of operations in
$\{\omega_{A}^{j,i}\}$; let
\[
z=\alpha_{A}\left(  C_{m+n-4}(\partial KK_{n,m})\right)
\]

\item[2.] Define $\beta$ on each monomial in $C_{\ast}\left(  \partial
JJ_{n,m}\smallsetminus\text{int}KK_{n,m}\times1\right)  $ to be its
corresponding fraction product of operations and maps in $\{\omega_{A}%
^{j,i},g_{i}^{j},\omega_{B}^{j,i}\}$; let
\[
\varphi=\beta\left(  C_{m+n-3}\left(  \partial JJ_{n,m}\smallsetminus
\text{int}KK_{n,m}\times1\right)  \right)
\]

\item[3.] Then $\tilde{g}\left(  z\right)  =\nabla\varphi$ implies $\left[
z\right]  =0$; choose a cochain $b$ such that $\nabla b=z$\smallskip

\item[4.] Note that $\nabla\left(  \tilde{g}\left(  b\right)  -\varphi\right)
=\nabla\tilde{g}\left(  b\right)  -\tilde{g}\left(  z\right)  =\tilde
{g}\left(  \nabla b-z\right)  =0$; choose a cocycle
\[
u\in\tilde{g}_{\ast}^{-1}\left[  \tilde{g}\left(  b\right)  -\varphi\right]
\]

\item[5.] Define $\alpha_{A}(\theta_{m}^{n})=\omega_{A}^{n,m}:=b-u$\smallskip

\item[6.] Define $\beta\left(  \left(  \mathfrak{f}_{1}^{1}\right)  ^{\otimes
n}\theta_{m}^{n}\right)  =g^{\otimes n}\omega_{A}^{n,m}$\smallskip

\item[7.] Note that $\left[  \tilde{g}\left(  \omega_{A}^{n,m}\right)
-\varphi\right]  =\left[  \tilde{g}\left(  b-u\right)  -\varphi\right]
=\left[  \tilde{g}\left(  b\right)  -\varphi\right]  -\left[  \tilde{g}\left(
u\right)  \right]  =0$; choose a cochain $g_{m}^{n}$ such that
\[
\nabla g_{m}^{n}=\tilde{g}\left(  \omega_{A}^{n,m}\right)  -\varphi=g^{\otimes
n}\omega_{A}^{n,m}-\varphi
\]

\item[8.] Define $\beta\left(  \mathfrak{f}_{m}^{n}\right)  =g_{m}^{n}$
\end{enumerate}

This completes the induction.$\medskip$

Let us apply the Transfer Algorithm to compute an induced $A_{\infty}$-algebra
structure on the cohomology of a DGA, which cannot be computed using the
classical BPL. The DGA $B$ considered here has no Hodge decomposition, its
homology $H=H^{\ast}\left(  B\right)  $ is not free, and the given homology
isomorphism $g:H\rightarrow B$ has no right-homotopy inverse.

Consider the tensor algebra $T^{a}M$ of the $\mathbb{Z}$-cochain complex
\[%
\begin{array}
[c]{cccccccccccc}%
\left(  M,d\right)  \text{\ } : & \mathbb{Z} & \rightarrow & 0 & \rightarrow &
\mathbb{Z}_{2}\oplus\mathbb{Z}_{2} & \rightarrow & \mathbb{Z}_{4} &
\rightarrow & \mathbb{Z}_{2} & \rightarrow & 0\\
& 1\mathbb{\!} &  & \mathbb{\!} &  & \left(  a,b\right)  \mathbb{\!} &
\mathbb{\!} & c & \mathbb{\!\mapsto} & x &  & \\
&  &  &  &  & b & \mathbb{\!\mapsto\!} & 2c\mathbb{\!} &  &  &  &
\end{array}
\]
and form the quotient $B=T^{a}M\left/  \left(  a^{2}+x,xc+cx,(ac+ca)^{2}%
,c^{2},tb,bt\right)  ,\right.  $ $\left\vert t\right\vert >0.$ Although the
DGA $B$ is not commutative, we do have $a\left(  ac+ca\right)  =\left(
ac+ca\right)  a$ and $\left(  ac\right)  ^{2}=\left(  ca\right)  ^{2}.$ Note
that $B$ has no Hodge decomposition since $c$ is not a cocycle, $2c$ is a
coboundary, and $\mathbb{Z}_{4}$ does not split as $\mathbb{Z}_{2}%
\oplus\mathbb{Z}_{2}$. Furthermore,
\[
H^{n}\left(  B\right)  =\left\{
\begin{array}
[c]{ll}%
\mathbb{Z} & n=0,\\
\mathbb{Z}_{2} & n=2,5,7\\
0 & \text{otherwise}%
\end{array}
\right.
\]
and $H=H^{\ast}\left(  B\right)  $ is not free. Define $g:H\rightarrow B$ by
$g\left(  1\right)  =1$ and
\[%
\begin{array}
[c]{lll}%
u=\left[  a\right]  & \mapsto & a\\
v=[ac+ca] & \mapsto & ac+ca\\
w=\left[  a\left(  ac+ca\right)  \right]  & \mapsto & a\left(  ac+ca\right)  .
\end{array}
\]
Then $g$ has no right-homotopy inverse $f$ since $gf\left(  b\right)  =\lambda
a$ implies $b-\lambda a=\left(  1-gf\right)  \left(  b\right)  =\left(
sd+ds\right)  \left(  b\right)  =sd\left(  b\right)  =s\left(  2c\right)
=2s\left(  c\right)  =0,$ which is a contradiction. To compute the induced
multiplication $\mu_{H}$, consider the following bases for $H$ and $H\otimes
H:$%
\[%
\begin{tabular}
[c]{|c|c|c|c|c|c|c|c|c|}\hline
& $2$ & $4$ & $5$ & $7$ & $9$ & $10$ & $12$ & $14$\\\hline
$H$ & $u$ &  & $v$ & $w$ &  &  &  & \\\hline
$H\otimes H$ &  & $u|u$ &  & $u|v,$ $v|u$ & $u|w,$ $w|u$ & $v|v$ & $v|w,$
$w|v$ & $w|w$\\\hline
\end{tabular}
\
\]
Ignoring the unit $1$ and evaluating $\tilde{g}$ on the basis $\left\{
w\partial_{u|v},w\partial_{v|u}\right\}  $ for $Hom^{0}\left(  H^{\otimes
2},H\right)  $ we have%
\[
\tilde{g}\left(  w\left(  \partial_{u|v}+\partial_{v|u}\right)  \right)
=g\left(  w\right)  \left(  \partial_{u|v}+\partial_{u|v}\right)  =\mu\left(
g\otimes g\right)  .
\]
Now thinking of $w\left(  \partial_{u|v}+\partial_{v|u}\right)  $ as a class
in $H^{\ast}\left(  Hom^{\ast}\left(  H^{\otimes2},H\right)  \right)  $ we
have%
\[
\tilde{g}_{\ast}\left[  w\left(  \partial_{u|v}+\partial_{v|u}\right)
\right]  =\left[  g\left(  w\right)  \left(  \partial_{u|v}+\partial
_{u|v}\right)  \right]  =\left[  \mu\left(  g\otimes g\right)  \right]  .
\]
Define $\mu_{H}=\left(  \tilde{g}_{\ast}\right)  ^{-1}\left[  \mu\left(
g\otimes g\right)  \right]  =w\left(  \partial_{u|v}+\partial_{v|u}\right)  ;$
then $uv=vu=w$ and $\mu_{H}$ is associative. To extend $g$ to an $A\left(
2\right)  $-map, let $\mu$ denote the multiplication in $B$ and consider the
expression
\begin{align*}
z  &  =\mu\left(  g\otimes g\right)  -g\mu_{H}=a^{2}\partial_{u|u}+\left(
a^{3}c+ca^{3}\right)  \left(  \partial_{u|w}+\partial_{w|u}\right) \\
&  =dc\partial_{u|u}+d\left(  cac\right)  \left(  \partial_{u|w}%
+\partial_{w|u}\right)  .
\end{align*}
Then $\nabla\left(  c\partial_{u|u}+cac\left(  \partial_{u|w}+\partial
_{w|u}\right)  \right)  =z$ and we define $g_{2}=c\partial_{u|u}+cac\left(
\partial_{u|w}+\partial_{w|u}\right)  $ so that
\[
\nabla g_{2}=\mu\left(  g\otimes g\right)  -g\mu_{H}.
\]
Thus $g$ is homotopy multiplicative. To compute the induced associator
$\mu_{H}^{3},$ consider the following bases for $H$ and $H^{\otimes3}:$%
\[%
\begin{tabular}
[c]{|c|c|c|c|c|c|}\hline
& $2$ & $5$ & $6$ & $7$ & $\cdots$\\\hline
$H$ & $u$ & $v$ &  & $w$ & \\\hline
$H\otimes H\otimes H$ &  &  & $u|u|u$ &  & $\cdots$\\\hline
\end{tabular}
\ \ \ \
\]
Since $B$ has trivial higher order structure, we consider the cochain
\[
\varphi=\mu\left(  g_{2}\otimes g-g\otimes g_{2}\right)  +g_{2}\left(  \mu
_{H}\otimes\mathbf{1}-\mathbf{1}\otimes\mu_{H}\right)  ,
\]
which vanishes on $H^{\otimes3}$ except $\varphi\left(  u|u|u\right)
=ac+ca=g\left(  v\right)  ;$ thus $\varphi=g\left(  v\right)  \partial
_{u|u|u}.$ Since $z=\mu_{H}\left(  \mu_{H}\otimes\mathbf{1}-\mathbf{1}%
\otimes\mu_{H}\right)  =0 $, every cocycle $b \in Hom^{-1}\left(  H^{\otimes
3},H\right)  $ satisfies $\nabla b=z$. Since $v\partial_{u|u|u}$ is the only
candidate, we set $b=v\partial_{u|u|u};$ then
\[
\left[  \tilde{g}\left(  b\right)  -\varphi\right]  =\left[  g\left(
v\right)  \partial_{u|u|u}-\varphi\right]  =\left[  0\right]  .
\]
Choose $u=0\in\tilde{g}_{\ast}^{-1}\left[  \tilde{g}\left(  b\right)
-\varphi\right]  $ and define $\mu_{H}^{3}=v\partial_{u|u|u};$ then $\mu
_{H}^{3}\left(  u|u|u\right)  =v.$ Finally, since $\varphi-g\mu_{H}^{3}%
\equiv0; $ we may set $g_{n}=0$ and $\mu_{H}^{n}=0$ for all $n\geq4$ to obtain
an induced $A_{\infty}$-algebra structure $\left(  H,\mu_{H},\mu_{H}%
^{3}\right)  $ and a map $G=g+g_{2}$ of $A_{\infty}$-algebras.

\section{A Topological Example}

Let $\mathbf{k}$ be a field. Given a space $X,$ let $S_{\ast}\left(  \Omega
X;\mathbf{k}\right)  $ denote the singular chains on the space of (base
pointed) Moore loops on $X,$ and choose a homology isomorphism $g:H_{\ast
}\left(  \Omega X;\mathbf{k}\right)  \rightarrow S_{\ast}\left(  \Omega
X;\mathbf{k}\right)  .$ Since $H=H_{\ast}\left(  \Omega X;\mathbf{k}\right)  $
is free and $S=S_{\ast}\left(  \Omega X;\mathbf{k}\right)  $ is a Hopf
algebra, the induced map $\tilde{g}:\mathcal{E}nd_{TH}\rightarrow U_{H,S}%
$\textit{\ }is a homology isomorphism by Proposition 1, and the Transfer
Algorithm induces an $A_{\infty}$-bialgebra structure on $H.$ Let us apply
this fact to a particular space $X$ and identify a non-trivial operation
$\omega_{2}^{2}:H\otimes H\rightarrow H\otimes H.$

Given a $1\,$-connected DGA $\left(  A,d_{A}\right)  $ over $\mathbb{Z}_{2},$
the bar construction of $A$, denoted by $BA,$ is the cofree DGC $T^{c}\left(
\downarrow\overline{A}\right)  $ with differential $d\ $and coproduct $\Delta$
defined as follows:\ Let $\left[  x_{1}|\cdots|x_{n}\right]  $ denote the
element $\left.  \downarrow x_{1}\otimes\cdots\otimes\downarrow x_{n}\right.
\in BA;$ then
\[
d\left[  x_{1}|\cdots|x_{n}\right]  =\sum_{i=1}^{n}\left[  x_{1}|\cdots
|dx_{i}|\cdots|x_{n}\right]  +\sum_{i=1}^{n-1}\left[  x_{1}|\cdots
|x_{i}x_{i+1}|\cdots|x_{n}\right]  ;
\]
\
\[
\Delta\left[  x_{1}|\cdots|x_{n}\right]  =\left[  \ \right]  \otimes\left[
x_{1}|\cdots|x_{n}\right]  +\left[  x_{1}|\cdots|x_{n}\right]  \otimes\left[
\ \right]  +\sum_{i=1}^{n-1}\left[  x_{1}|\cdots|x_{i}\right]  \otimes\left[
x_{i}|\cdots|x_{n}\right]  .
\]

Consider the space $Y=\left(  S^{2}\times S^{3}\right)  \vee\Sigma
\mathbb{C}P^{2},$ multiplicative generators $\bar{a}_{i}\in H^{i}\left(
S^{i};\mathbb{Z}_{2}\right)  ,$ $\bar{b}\in H^{3}\left(  \Sigma\mathbb{C}%
P^{2};\mathbb{Z}_{2}\right)  $ and $Sq^{2}\bar{b}\in H^{5}\left(
\Sigma\mathbb{C}P^{2};\mathbb{Z}_{2}\right)  $, a map $f:Y\rightarrow
K(\mathbb{Z}_{2},5)$ such that $f^{\ast}(\iota_{5})=\bar{a}_{2}\bar{a}%
_{3}+Sq^{2}\bar{b},$ and the pullback $p:X\rightarrow Y$ of the following path
fibration:
\[%
\begin{array}
[c]{ccccc}%
K\left(  \mathbb{Z}_{2},4\right)  & \longrightarrow & X\medskip &
\longrightarrow & \mathcal{L}K\left(  \mathbb{Z}_{2},5\right) \\
&  & p\downarrow\medskip\text{ } &  & \downarrow\\
&  & Y & \overset{f}{\longrightarrow} & K\left(  \mathbb{Z}_{2},5\right)  .\\
&  & \bar{a}_{2}\bar{a}_{3}+Sq^{2}\bar{b} & \underset{f^{\ast}}{\longleftarrow
} & \iota_{5}%
\end{array}
\]
Let $a_{i}=p^{\ast}(\bar{a}_{i})$ and $b=p^{\ast}(\bar{b})$; then
$A\!=\!H^{\ast}\!\left(  X;\mathbb{Z}_{2}\right)  \!=\!\left\{  1,a_{2}%
,a_{3},b,a_{2}a_{3}=Sq^{2}b,\ldots\right\}  .$

Form the bar construction $BA;$ since $H=H^{\ast}\left(  BA\right)  \approx
H^{\ast}\left(  \Omega X;\mathbb{Z}_{2}\right)  $ as coalgebras, $\left(
BA,d,\Delta\right)  $ is a DG coalgebra model for cochains on $\Omega X$. In
\cite{Baues}, H.-J. Baues identified a compatible multiplication
$\mu:BA\otimes BA\rightarrow BA$ and a DG\ Hopf algebra model $\left(
BA,d,\Delta,\mu\right)  $ in the following way: The twisting in $X$ induces
Steenrod's $\smile_{1}:A\otimes A\rightarrow A,$ which acts non-trivially via
$b\smile_{1}b=a_{2}a_{3}$ and the induced map $\phi:BA\otimes BA\rightarrow A$
acts non-trivially via%
\[
\phi\left(  \left[  x\right]  \otimes\left[  \ \right]  \right)  =\phi\left(
\left[  \ \right]  \otimes\left[  x\right]  \right)  =x\text{ and }\phi\left(
\left[  b\right]  \otimes\left[  b\right]  \right)  =b\smile_{1}b=a_{2}a_{3}.
\]
(cf. \cite{Voronov}, \cite{KS1}). Consider the tensor product of coalgebras
$BA\otimes BA$ with coproduct $\psi=\sigma_{2,2}\left(  \Delta\otimes
\Delta\right)  $ and define
\[
\mu:=\sum\nolimits_{k\geq0}\text{ }\underset{k+1\text{ factors}%
}{\underbrace{\left(  \downarrow\phi\otimes\cdots\otimes\downarrow\phi\right)
_{\mathstrut}}}\text{ }\bar{\psi}^{\left(  k\right)  },
\]
where $\bar{\psi}^{(0)}=\mathbf{1}$, $\bar{\psi}^{\left(  k\right)  }=\left(
\bar{\psi}\otimes\mathbf{1}^{\otimes k-1}\right)  \cdots\left(  \bar{\psi
}\otimes\mathbf{1}\right)  \bar{\psi}$ for $k>0,$ and $\bar{\psi}$ is the
reduced coproduct. Then for example, $\mu\left(  \lbrack b]\otimes\left[
b\right]  \right)  =\left[  a_{2}a_{3}\right]  .$

Let $\mu_{H}$ be the multiplication on $H$ induced by $\mu$ and consider the
classes $\alpha_{i}=\text{cls}\left[  a_{i}\right]  ,\beta=\text{cls}\left[
b\right]  \in H.$ Choose a cocycle-selecting map $g:H\rightarrow BA$ such that
$g\left(  \text{cls}\left[  x_{1}|\cdots|x_{n}\right]  \right)  =\left[
x_{1}|\cdots|x_{n}\right]  .$ Then $\mu\left(  \left[  b\right]
\otimes\left[  b\right]  \right)  =\left[  a_{2}a_{3}\right]  =d\left[
a_{2}|a_{3}\right]  $ implies $\mu_{H}\left(  \beta\otimes\beta\right)  =0$
and $\left(  g\mu_{H}+\mu\left(  g\otimes g\right)  \right)  \left(
\beta\otimes\beta\right)  =\left[  a_{2}a_{3}\right]  $ . Nevertheless, by the
Transfer Theorem, there is a cochain homotopy $g_{2}^{1}:H\otimes H\rightarrow
BA$ satisfying the relation $\nabla g_{2}^{1}=g\mu_{H}+\mu\left(  g\otimes
g\right)  $ such that $g_{2}^{1}\left(  \beta\otimes\beta\right)  =\left[
a_{i}|a_{5-i}\right]  $ for some $i\in\left\{  2,3\right\}  ;$ and in
particular, we may choose
\[
g_{2}^{1}\left(  \beta\otimes\beta\right)  =\left[  a_{2}|a_{3}\right]
\]
since either choice gives rise to isomorphic structures. Let $\Delta_{H}$ be
the coproduct induced by $\Delta;$ then $\left\{  \Delta g+\left(  g\otimes
g\right)  \Delta_{H}\right\}  \left(  \beta\right)  =0 $ since $\beta$ is
primitive. By the Transfer Theorem, there is a cochain homotopy $g_{1}%
^{2}:H\rightarrow BA\otimes BA$ satisfying the relation $\nabla g_{1}%
^{2}=\Delta g+\left(  g\otimes g\right)  \Delta_{H}$ such that $\nabla
g_{1}^{2}\left(  \beta\right)  =0.$ Thus $g_{1}^{2}\left(  \beta\right)
=\lambda\otimes\left[  a_{2}\right]  +\left[  a_{2}\right]  \otimes\rho$ \ for
some $\lambda,\rho\in\mathbb{Z}_{2};$ and in particular, we may choose%
\[
g_{1}^{2}\left(  \beta\right)  =0.
\]
Again by the Transfer Theorem, there is a cochain homotopy $g_{2}^{2}:H\otimes
H \rightarrow BA \otimes BA$ satisfying the following relation on $JJ_{2,2}:$
\begin{align}
\nabla g_{2}^{2}  &  =\left(  \mu\otimes\mu\right)  \sigma_{2,2}\left(  \Delta
g\otimes g_{1}^{2}+g_{1}^{2}\otimes\left(  g\otimes g\right)  \Delta
_{H}\right) \label{one}\\
&  +\left(  \mu\left(  g\otimes g\right)  \otimes g_{2}^{1}+g_{2}^{1}\otimes
g\mu_{H}\right)  \sigma_{2,2}\left(  \Delta_{H}\otimes\Delta_{H}\right)
\nonumber\\
&  +\omega_{BA}^{2,2}\left(  g\otimes g\right)  +\left(  g\otimes g\right)
\omega_{H}^{2,2}+\Delta g_{2}^{1}+g_{1}^{2}\mu_{H}.\nonumber
\end{align}
The component $\omega_{BA}^{2,2}\left(  g\otimes g\right)  $ vanishes since
$BA$ has trivial higher order structure; the non-triviality of $\left(
g\otimes g\right)  \omega_{H}^{2,2}$ is to be determined.

Let us evaluate relation (\ref{one}) at $\beta\otimes\beta.$ First, $g_{1}%
^{2}\mu_{H}\left(  \beta\otimes\beta\right)  =0$ by the observation above, and
$\left(  \mu\otimes\mu\right)  \sigma_{2,2}\left(  \Delta g\otimes g_{1}%
^{2}+g_{1}^{2}\otimes\left(  g\otimes g\right)  \Delta_{H}\right)  \left(
\beta\otimes\beta\right)  =0$ by our choice of $g_{1}^{2}$. Second, $\left(
\mu\left(  g\otimes g\right)  \otimes g_{2}^{1}+g_{2}^{1}\otimes g\mu
_{H}\right)  \sigma_{2,2}\left(  \Delta_{H}\otimes\Delta_{H}\right)  \left(
\beta\otimes\beta\right)  =\left[  \ \right]  \otimes g_{2}^{1}\left(
\beta\otimes\beta\right)  +g_{2}^{1}\left(  \beta\otimes\beta\right)
\otimes\left[  \ \right]  =\left(  \Delta+\overline{\Delta}\right)  g_{2}%
^{1}\left(  \beta\otimes\beta\right)  .$ Thus relation (\ref{one}) reduces to
$\nabla g_{2}^{2}(\beta\otimes\beta)=\left(  g\otimes g\right)  \omega
_{H}^{2,2}\left(  \beta\otimes\beta\right)  +\overline{\Delta}g_{2}^{1}\left(
\beta\otimes\beta\right)  =\left(  g\otimes g\right)  \omega_{H}^{2,2}\left(
\beta\otimes\beta\right)  +\left[  a_{2}\right]  \otimes\left[  a_{3}\right]
,$ and we conclude that
\[
\omega_{H}^{2,2}\left(  \beta\otimes\beta\right)  =\alpha_{2}\otimes\alpha
_{3}.
\]
Thus the topology of the total space $X$ in the fibration $p:X\rightarrow
\left(  S^{2}\times S^{3}\right)  \vee\Sigma\mathbb{C}P^{2}$ above induces a
nontrivial $2$-in/$2$-out operation $\omega_{H}^{2,2}$ on $H=H^{\ast}\left(
\Omega X;\mathbb{Z}_{2}\right)  .$ A variation of this example, with a
nontrivial topologically induced $2$-in/$n$-out operation on loop cohomology
for each $n\geq2$, appears in \cite{SU5}.\vspace{0.2in}

\textbf{Acknowledgement. }I wish to thank Samson Saneblidze for sharing his
deep insights, which enabled this work, and Jim Stasheff for his thoughtful
questions and stimulating discussions.

\end{document}